\documentclass[12pt,a4paper]{article}

\usepackage{amssymb,amsmath,amsfonts,amsthm,enumerate}
\usepackage{graphicx}

\usepackage[utf8]{inputenc}
\usepackage[english]{babel}

\numberwithin{equation}{section}
\newtheorem{theorem}{Theorem}[section]

\theoremstyle{definition}

\newtheorem{definition}[theorem]{Definition}
\newtheorem{remark}[theorem]{Remark}

\DeclareMathOperator*{\esup}{ess\,sup}

\title{Weak solvability of a boundary value problem\\ for a parabolic equation with a global-in-time term
\\ that contains a weighted integral}
\author{V.\,N.\,Starovoitov\\
\small Lavrentyev Institute of Hydrodynamics, Novosibirsk, Russian Federation\\
\small E-mail:\; starovoitov@hydro.nsc.ru}
\date{}

\begin{document}

\maketitle

\abstract{
This paper deals with a parabolic partial differential
equation that includes a non-linear nonlocal in time term. This term
is the product of a so-called interaction potential and the solution of the problem.
The interaction potential depends on a weighted integral of the solution over the entire time interval, where the
problem is considered, and satisfies fairly general conditions.
Namely, it is assumed to be a continuous bounded from below function that can behave arbitrarily at infinity.
This fact implies that the interaction term is not a lower order term in the equation.
The weak solvability of the initial boundary value problem for this equation is proven.
The proof does not use any continuity properties of the solution with respect to time and is based on the energy
estimate only.
}

\bigskip\noindent
\textbf{Key words:} nonlocal in time parabolic equation, weighted integral, initial boundary value problem, solvability

\bigskip\noindent
\textbf{2010 Mathematics Subject Classification:} 35K58, 35Q92

\section{Introduction}\label{sec1}
Let $\varOmega$ be a bounded domain in $\mathbb{R}^n$, $n\ge 2$, with a Lipschitz boundary $\partial\varOmega$.
In the space-time cylinder $\varOmega_T=\varOmega\times (0,T)$, $T\in(0,\infty)$, we consider the following differential equation:
\begin{equation}\label{1.1}
\partial_t u -\Delta u +\varphi\Big(\int_0^T \alpha(\cdot,s)\, u(\cdot,s)\,ds\Big)\, u =f,
\end{equation}
where $u=u(x,t)$ is an unknown scalar function, $x=(x_1,\ldots,x_n)$ the vector of the spatial variables
in $\mathbb{R}^n$, $t$ the time variable on $[0,T]$, $\varphi: \mathbb{R}\to \mathbb{R}$ a
function that will be specified below, $f:\varOmega_T\to \mathbb{R}$ a prescribed function.  We suppose that
the following boundary and initial conditions are satisfied:
\begin{equation}\label{1.2}
u(x,t)=0\quad \text{for}\quad x\in\partial\varOmega,\quad t\in [0,T],
\end{equation}
\begin{equation}\label{1.3}
u(x,0)=u_0(x)\quad \text{for}\quad x\in\varOmega,
\end{equation}
where the function $u_0:\varOmega\to \mathbb{R}$ is prescribed.

Equation \eqref{1.1} contains a non-local in time $t$ term with the integral over the
whole interval $(0, T)$, where the problem is considered. For this reason, we call this equation global-in-time.
A problem of this type arose when modeling the chaotic dynamics of a polymer molecule (polymer chain) in an aqueous solution
\cite{StSt}. The function $u(x,t)$ corresponds to the density of probability that the
$t$-th segment of the chain is in a certain region of the space. The role of time in the equation
is played by the arc length parameter along the chain whose total length is $T$.
The equation contains a term responsible for the interaction of chain's segments and the function
$\varphi$ is called the \emph{interaction potential}. Since each segment interacts with all others
through the surrounding fluid, $\varphi$ depends on the integral of the density of probability over the entire chain,
i.e., over the entire time interval $(0,T)$.

Equations with a global-in-time term appear also in other fields, for example, in the filtration theory \cite{L},
in the population dynamics (see \cite{Webb,W1} and the references therein).
Notice that the equation in the population dynamics is ultra-parabolic
and the age of individuals plays the role of the second time. In the steady case, the equation becomes parabolic.
The global terms contain the integral with respect to the age and can be in the data of the problem and in the equation as well.
There are a lot of works, where a global-in-time term stands in the data of the problem.
These are problems with non-local boundary or initial conditions. The latter means that the initial data depend on
the time integral of the solution. This topic is widely represented in the literature, and there is no problem to find
relevant publications. The problems with global-in-time data are very different from ours.

The problem we study looks like parabolic one, however it has some unusual features.
First of all, since  the potential $\varphi$ in equation \eqref{1.1} is global-in-time, the state of the system
depends not only on the past but also on the future. That is, the causality principle is violated. Besides, from a mathematical point
of view, the solution of a nonlinear parabolic problem is commonly being constructed locally in time and is extended
afterwards. In our case, this procedure is impossible. Finally, as a rule, the local-in-time uniqueness of the solution implies
the global one. We cannot prove the uniqueness without restrictions on $T$.

There are several papers devoted to problems close to \eqref{1.1}--\eqref{1.3}. In \cite{St1,St2},
the weak solvability of problem \eqref{1.1}--\eqref{1.3} with $\alpha\equiv 1$ and $f=0$ is proven. Notice that
the smoothness of the obtained in this papers weak solution can be easily improved for better initial data.
The strong solvability of the same problem was proven in \cite{W2}, where the semigroup approach was employed.
In \cite{D}, the authors consider the problem with the Laplace operator replaced by a more complex non-local
operator that can be thought as the fractional Laplacian. The uniqueness of the solution is proven for sufficiently
small $T$ in all these works. The smallness condition means that there is a restriction for $T$ which depends on the
data of the problem.

In this paper, we consider a more general problem compared with the previous works. At first, the global-in-time term
contains a weighted integral and the smoothness of the weight $\alpha$ is assumed to be not sufficient to apply the
technique from \cite{St1,St2}. Secondly, equation \eqref{1.1} contains the right-hand side $f$.
The point is that the proof of the main results of the paper employs only the energy estimate and does not use any
smoothness properties of the solution with respect to time. The right-hand side $f$ prevents or makes it difficult to
obtain good time-estimates. Finally, the potential $\varphi$ is assumed to be just a continuous bounded from below function.
There are no conditions on the growth of this function at infinity, no convexity, monotonicity, or differentiability conditions.
In \cite{St2},  the assumption on the potential $\varphi$ also does not contain the growth conditions, but $\varphi$
satisfies other requirements that made it possible to obtain a higher degree of integrability of $\varphi$ than in the presented paper.
The absence of the growth conditions at infinity leads to the fact that the interaction term is not a lower order term
in equation \eqref{1.1}. Notice that in \cite{W2} this term has a lower order.

Although the continuous function $\varphi$ is assumed to be just bounded from below, by standard arguments
the problem can be reduced to the problem with a positive potential. Really,
suppose that $\varphi(\xi)\ge -K$ for some positive number $K$ and all $\xi\in \mathbb{R}$.
If $u$ is a solution of problem~\eqref{1.1}--\eqref{1.3}, then $\bar{u}(x,t)=e^{-Kt}\,u(x,t)$ is a solution
of the following problem
$$
\partial_t \bar{u} -\Delta \bar{u} +\bar{\varphi}(\bar{v})\, \bar{u} =\bar{f},
\quad \bar{u}|_{\partial\varOmega}=0,\quad \bar{u}|_{t=0}=u_0,
$$
where
$$
\bar{\varphi}=\varphi +K,\quad
\bar{v}(x)=\int_0^T \bar{\alpha}(x,t)\, \bar{u}(x,t)\,dt,
\quad \bar{\alpha}(x,t)=e^{Kt}\,\alpha(x,t), \quad \bar{f}(x,t)=e^{-Kt}\,f(x,t).
$$
Thus, we have the same problem but with the non-negative potential $\bar{\varphi}$.
In what follows, we will assume that $\varphi$ is a non-negative function.

We will use the standard Lebesgue and Sobolev spaces $L^p(\varOmega)$, $H^1_0(\varOmega)$,
$L^q(0,T;L^p(\varOmega))$, and $L^q(0,T;H^1_0(\varOmega))$, where $p,q\in [1,\infty]$.
As usual, $H^{-1}(\varOmega)$ is the dual space of $H^1_0(\varOmega)$
with respect to the pivot space $L^2(\varOmega)$. The norm in $L^2(\varOmega)$ will be denoted by $\|\cdot\|$.

\begin{definition}\label{t1.1}
Let $\varphi: \mathbb{R}\to \mathbb{R}$ be a non-negative continuous function, $T\in (0,\infty)$,
$\alpha\in L^1(0,T; L^2(\varOmega))$, $f\in L^1(0,T; L^2(\varOmega))$, and $u_0\in L^2(\varOmega)$.
A function $u:\varOmega_T\to \mathbb{R}$ is said to be a \emph{weak solution of problem \eqref{1.1}--\eqref{1.3}}, if
\begin{enumerate}
\item
$u\in L^2(0,T;H^1_0(\varOmega))$ and $\varphi(\zeta)\,u\in L^1(\varOmega_T)$, where
$\zeta=\int_0^T \alpha (\cdot,t)\, u(\cdot,t)\,dt$;
\item
the following integral identity
$$
\int_0^T \int_\varOmega \big(u\,\partial_t h - \nabla u\cdot\nabla h - \varphi(v)\,u\, h +f\,h\big)\,dx\,dt +
\int_\varOmega u_0 h_0\, dx=0
$$
holds for an arbitrary smooth in the closure of $\varOmega_T$ function $h$ such that $h(x,t)=0$
for $x\in\partial\varOmega$ and for $t=T$. Here, $h_0=h|_{t=0}$.
\hfill\textbullet
\end{enumerate}
\end{definition}

The main result of the paper is Theorem~\ref{t3.1} that states the existence of the weak solution of
problem~\eqref{1.1}--\eqref{1.3}. In the next section, we consider the case where the potential is
assumed to be bounded and prove Theorem~\ref{t2.1} that is an auxiliary result for Theorem~\ref{t3.1}.
We emphasize once again that in the proofs of these theorems we do not use any continuity properties of the solution in time.
This leads to certain difficulties, for example, in proving the energy estimate \eqref{3.1} (see Step~4 of the proof of
Theorem~\ref{t3.1}) which is usually obtained automatically from
the corresponding estimate for the approximate solutions by using the weak semicontinuity of norms in Banach spaces.

\section{Problem with a bounded potential}\label{sec2}
In this section, we consider the case of the problem where the potential $\varphi:\mathbb{R}\to \mathbb{R}$
is a continuous bounded function.
\begin{theorem}\label{t2.1}
Assume that $T\in (0,\infty)$,  $\alpha\in L^1(0,T; L^2(\varOmega))$,
$f\in L^1(0,T; L^2(\varOmega))$, $u_0\in L^2(\varOmega)$, and $\varphi:\mathbb{R}\to \mathbb{R}$
is a continuous function such that
$0\le \varphi(\xi)\le K$ for all $\xi\in \mathbb{R}$, where $K$ is a positive number.
Then there exists a weak solution $u\in L^\infty(0,T;L^2(\varOmega))\cap L^2(0,T;H^1_0(\varOmega))$
of problem~\eqref{1.1}--\eqref{1.3}.
\end{theorem}
\noindent
\textsc{Proof.} To prove the solvability of the problem, we employ the Schauder fixed point theorem. Let $B_R$ be the
ball of the radius $R$  centered at the zero in the function space $L^1(\varOmega)$. The number $R$ will be
determined later. Let us define the mapping $\varPsi:B_R\to B_R$ that is referred to in the Schauder theorem.
For every $w\in B_R$, we define $u_w$ as the solution of the following problem:
\begin{equation}\label{2.1}
\partial_t u_w -\Delta u_w +\varphi(w)\, u_w =f,
\quad u_w|_{\partial\varOmega}=0,\quad u_w|_{t=0}=u_0.
\end{equation}
As it follows from the classical theory of parabolic equations (see, e.g., \cite[Ch.\,7]{Evans}, this problem has a unique weak solution.
Now, we set
$$
\varPsi(w)=\int_0^T \alpha(\cdot,t)\, u_w(\cdot,t)\,dt.
$$
It necessary to prove that there exists $R>0$ such that $\varPsi(B_R)\subset B_R$ and $\varPsi:B_R\to L^1(\varOmega)$
is a compact continuous map.

The usual energy estimate for problem \eqref{2.1} looks as follows:
$$
\|u_w(\cdot,t)\|^2+\int_0^t\|\nabla u_w(\cdot,s)\|^2\,ds \le
\Big(\|u_0\| +\int_0^t \|f(\cdot,s)\|\,ds\Big)^2+\frac{1}{2}\,\Big(\int_0^t \|f(\cdot,s)\|\,ds\Big)^2
$$
for almost all $t\in [0,T]$. We have used the positiveness of $\varphi$ although  a similar estimate
holds for $\varphi$ which is just bounded. Thus,
\begin{equation}\label{2.2}
\esup_{t\in[0,T]}\|u_w(\cdot,t)\|^2+\int_0^T\|\nabla u_w(\cdot,s)\|^2\,ds \le C_1,
\end{equation}
where $C_1=\big(\|u_0\| +\int_0^T \|f(\cdot,s)\|\,ds\big)^2+\frac{1}{2}\,\big(\int_0^T \|f(\cdot,s)\|\,ds\big)^2$.
As a consequence of this inequality, we find that
$$
\|\varPsi(w)\|_1\le C_2\sqrt{C_1},
$$
where $C_2=\int_0^T \|\alpha(\cdot,s)\|\,ds$. Therefore, $\varPsi(B_R)\subset B_R$ with $R=C_2\sqrt{C_1}$.

In order to prove the continuity of the mapping $\varPsi$, let us take an arbitrary sequence $\{w_k\}$ in $B_R$ that converges
in $L^1(\varOmega)$ to some $w\in B_R$. Then $w_k\to w$ in measure on $\varOmega$ and, due to the Lebesgue
dominated convergence theorem, $\varphi(w_k)\to \varphi(w)$ in $L^p(\varOmega)$ for all $p\in [1,\infty)$.
The difference $u_{w_k}-u_w$ is a weak solution of the following problem:
\begin{gather*}
\partial_t (u_{w_k}-u_w) -\Delta (u_{w_k}-u_w) +\varphi(w)\, (u_{w_k}-u_w) +
\big(\varphi(w_k)-\varphi(w)\big)\, u_{w_k}=0,
\\
\quad (u_{w_k}-u_w)|_{\partial\varOmega}=0,\quad (u_{w_k}-u_w)|_{t=0}=0.
\end{gather*}
Therefore, it satisfies estimate \eqref{2.2} with $f=\big(\varphi(w_k)-\varphi(w)\big)\, u_{w_k}$ and $u_0=0$. Thus,
$$
\esup_{t\in[0,T]}\|(u_{w_k}-u_w)(\cdot,t)\|^2+\int_0^T\|\nabla (u_{w_k}-u_w)(\cdot,s)\|^2\,ds \le c_k,
$$
where $c_k=\frac{3}{2}\,\big(\int_0^T \|\big(\varphi(w_k)-\varphi(w)\big)\, u_{w_k}(\cdot,s)\|\,ds\big)^2$.
Due to the Sobolev embedding theorems, estimate \eqref{2.2} implies that
$$
\int_0^T \|u_{w_k}(\cdot,s)\|_{L^q(\varOmega)}^2\,ds\le C_3,
$$
for all $k\in \mathbb{N}$ and some independent of $k$ constant $C_3$, where $q=2n/(n-2)$,
if $n>2$, and $q=\infty$ otherwise. As a consequence of the H\"older inequality, we have
$$
\int_0^T \|\big(\varphi(w_k)-\varphi(w)\big)\, u_{w_k}(\cdot,s)\|\,ds
\le \|(\varphi(w_k)-\varphi(w)\|_{L^p(\varOmega)}\, \int_0^T \|u_{w_k}(\cdot,s)\|_{L^q(\varOmega)}\,ds,
$$
where $p=2q/(q-2)$. Thus, $c_k\to 0$ and $\esup_{t\in[0,T]}\|(u_{w_k}-u_w)(\cdot,t)\|\to 0$ as $k\to \infty$.
This means that $\varPsi(w_k)\to \varPsi(w)$ in $L^1(\varOmega)$ as $k\to\infty$. Really,
\begin{multline*}
\|\varPsi(w_k)- \varPsi(w)\|_{L^1(\varOmega)}\le \int_0^T \|\alpha(\cdot,s)\|\, \|(u_{w_k}-u_w)(\cdot,s)\|\,ds
\\
\le \int_0^T \|\alpha(\cdot,s)\|\,ds\, \esup_{t\in[0,T]}\|(u_{w_k}-u_w)(\cdot,t)\|\to 0
\quad\text{as}\quad k\to\infty.
\end{multline*}
Thus, the mapping $\varPsi$ is continuous on $B_R$.

It remains to prove the compactness of $\varPsi$. It is possible to draw arguments based on the available
higher smoothness of the solution to problem \eqref{2.1}. In fact, $u_w$ satisfies a stronger estimate than
\eqref{2.2}. However, we will show that this estimate is sufficient to prove the compactness of $\varPsi$.
Let $\{w_k\}$ be a bounded sequence in $B_R$. We need to show that there exists a subsequence $\{w_{k'}\}$
such that $\varPsi(w_{k'})$ converges in $L^1(\varOmega)$. According to estimate \eqref{2.2},
the sequence $\{u_{w_k}\}$ is bounded in $L^\infty(0,T;L^2(\varOmega))\cap L^2(0,T;H^1_0(\varOmega))$.
Besides that,  $\|\nabla\int_0^T h \,u_{w_{k}}\,dt\|\le C_4$ for
an arbitrary smooth function $h:\varOmega_T\to \mathbb{R}$. The constant $C_4$ depends, of course, on $h$.
Therefore, there exists a function $v\in L^\infty(0,T;L^2(\varOmega))\cap L^2(0,T;H^1_0(\varOmega))$
and a subsequence $\{u_{w_{k'}}\}$  such that $u_{w_{k'}}\to v$ as $k'\to\infty$
weakly in $L^2(0,T;H^1_0(\varOmega))$ and $*$-weakly in $L^\infty(0,T;L^2(\varOmega))$, and
$$
\int_0^T h \,u_{w_{k'}}\,dt\to \int_0^T h \,v\,dt\quad\text{in}\quad L^2(\varOmega).
$$

Let us take an arbitrary $\varepsilon>0$. There exists a smooth function $\alpha_\varepsilon:\varOmega_T\to \mathbb{R}$ such that
$\|\alpha-\alpha_\varepsilon\|_{L^1(0,T;L^2(\varOmega))}<\varepsilon$. If we take  $h=\alpha_\varepsilon$ in the previous
relation, then we find that
\begin{multline*}
\limsup_{k'\to\infty}\Big\|\int_0^T\alpha\,(u_{w_{k'}}-v)\,dt\Big\|_{L^1(\varOmega)}\le
\limsup_{k'\to\infty}\Big\|\int_0^T(\alpha-\alpha_\varepsilon)\,(u_{w_{k'}}-v)\,dt\Big\|_{L^1(\varOmega)}
\\
+\limsup_{k'\to\infty}\Big\|\int_0^T\alpha_\varepsilon (u_{w_{k'}}-v)\,dt\Big\|_{L^1(\varOmega)}
\le C_5\,\|\alpha-\alpha_\varepsilon\|_{L^1(0,T;L^2(\varOmega))} < C_5\,\varepsilon,
\end{multline*}
where the constant $C_5$ is an upper bound for the norm of  $u_{w_{k'}}-v$ in $L^\infty(0,T;L^2(\varOmega))$.
Since $\varepsilon$ is arbitrary,
$$
\lim_{k'\to\infty}\Big\|\int_0^T\alpha\,(u_{w_{k'}}-v)\,dt\Big\|_{L^1(\varOmega)}=0.
$$
Thus, the sequence $\{\varPsi(w_{k'})\}$ converges in $L^1(\varOmega)$.

So, all the conditions of the Schauder theorem is fulfilled and there exists $w\in L^1(\varOmega)$ such that
$\varPsi(w)=w$. Then $u=u_w$ is a weak solution of problem~\eqref{1.1}--\eqref{1.3} and the theorem is proven.
\hfill$\square$

\begin{remark}\label{t2.2}
In Theorem~\ref{t2.1}, we assumed that $f\in L^1(0,T; L^2(\varOmega))$. Without changing the proof, we could suppose
that  $f\in L^2(0,T; H^{-1}(\varOmega))$.  Really, the proof is based on the energy estimate \eqref{2.2}.
This estimate is also true with $C_1=\|u_0\|^2 +c^2 \int_0^T \|f(\cdot,s)\|^2_{H^{-1}(\varOmega)}\,ds$,
where the constant $c$ is such that $\|u\|_{H^1_0(\varOmega)}\le c\,\|\nabla u\|$.
\hfill$\bullet$
\end{remark}

\section{Problem with an unbounded potential}\label{sec3}
In this section, we consider a more general case where the potential $\varphi$
is a non-negative continuous function which is not necessarily bounded.
\begin{theorem}\label{t3.1}
Assume that $T\in (0,\infty)$,  $\alpha\in L^\infty(\varOmega; L^2(0,T))$,
$f\in L^1(0,T; L^2(\varOmega))$, $u_0\in L^2(\varOmega)$, and $\varphi:\mathbb{R}\to \mathbb{R}$
is a non-negative continuous function. Then there exists a weak solution $u\in
L^\infty(0,T;L^2(\varOmega))\cap L^2(0,T;H^1_0(\varOmega))$ of problem~\eqref{1.1}--\eqref{1.3} such that
\begin{equation}\label{3.1}
\esup_{t\in[0,T]}\|u(\cdot,t)\|^2+\int_0^T\|\nabla u(\cdot,t)\|^2\,dt
+\int_0^T \int_\varOmega \varphi(\zeta(x))\, u^2(x,t)\, dx\, dt
\le C_1,
\end{equation}
\begin{equation}\label{3.2}
\text{$\varphi(\zeta)$, $\varphi(\zeta)\zeta$, and $\varphi(\zeta)\zeta^2$ are in $L^1(\varOmega)$},
\end{equation}
\begin{equation}\label{3.3}
\text{$\varphi(\zeta)u$ and $\varphi(\zeta)u^2$ are in $L^1(\varOmega_T)$},
\end{equation}
where $\zeta=\int_0^T \alpha (\cdot,t)\, u(\cdot,t)\,dt\in L^2(\varOmega)$ and
$C_1=\big(\|u_0\| +\int_0^T \|f(\cdot,s)\|\,ds\big)^2+\frac{1}{2}\,\big(\int_0^T \|f(\cdot,s)\|\,ds\big)^2$.
\end{theorem}
\noindent
\textsc{Proof.}
For convenience, we split the proof into several steps.

\medskip\noindent
\textsc{Step~1.}
For every $k\in \mathbb{N}$, we define the following function:
$$
\varphi_k(\xi)=\begin{cases}
\varphi(\xi), & \varphi(\xi)\le k,
\\
k, & \text{otherwise}.
\end{cases}
$$
Let us denote by $u_k$ the weak solution of problem~\eqref{1.1}--\eqref{1.3} with the potential $\varphi_k$.
Since $\varphi_k$ is bounded and $\alpha\in L^\infty(\varOmega; L^2(0,T))\subset L^1(0,T; L^2(\varOmega))$,
the existence of $u_k$ is proven in the previous section (Theorem~\ref{t2.1}).
For brevity, we introduce the notation:
$$
\zeta_k=\int_0^T \alpha (\cdot,t)\, u_k(\cdot,t)\,dt.
$$
According to Definition~\ref{t1.1}, $u_k$ satisfies the integral identity
\begin{equation}\label{3.4}
\int_0^T \int_\varOmega \big(u_k\,\partial_t h - \nabla u_k\cdot\nabla h - \varphi_k(\zeta_k)\,u_k\, h
+ f\, h\big)\,dx\,dt + \int_\varOmega u_0 \,h_0\, dx=0
\end{equation}
for an arbitrary smooth in the closure of $\varOmega_T$ function $h$ such that $h(x,t)=0$
for $x\in\partial\varOmega$ and for $t=T$. The energy estimate for $u_k$ looks as follows:
\begin{equation}\label{3.5}
\esup_{t\in[0,T]}\|u_k(\cdot,t)\|^2+\int_0^T\|\nabla u_k(\cdot,t)\|^2\,dt
+\int_0^T \int_\varOmega \varphi_k(\zeta_k(x))\, u_k^2(x,t)\, dx\, dt
\le C_1,
\end{equation}
where the constant $C_1$ is the same as in \eqref{2.2} and \eqref{3.1}. This estimate implies that the sequence $\{u_k\}$
has a subsequence (denoted again by $\{u_k\}$) which converges $*$-weakly in $L^\infty(0,T;L^2(\varOmega))$
and weakly in $L^2(0,T;H^1_0(\varOmega))$ to a function $u$. Our goal is to pass to the limit in
\eqref{3.4} as $k\to\infty$. Clearly,
$$
\lim_{k\to\infty}\int_0^T \int_\varOmega \big(u_k\,\partial_t h - \nabla u_k\cdot\nabla h\big)\,dx\,dt
=\int_0^T \int_\varOmega \big(u\,\partial_t h - \nabla u\cdot\nabla h\big)\,dx\,dt.
$$
The only problem is to prove that
\begin{equation}\label{3.6}
\lim_{k\to\infty}\int_0^T \int_\varOmega \varphi_k(\zeta_k)\,u_k\, h\,dx\,dt=
\int_0^T \int_\varOmega \varphi(\zeta)\,u\, h\,dx\,dt.
\end{equation}
Notice that we also need to justify that the right-hand side in the last relation makes sense.

\medskip\noindent
\textsc{Step~2.}
The reasoning at the end of the proof of Theorem~\ref{t2.1} implies that
$$
\zeta_k\to\zeta\quad\text{in $L^1(\varOmega)$ as $k\to\infty$.}
$$
Consequently, there exists a subsequence (denoted again by $\{\zeta_k\}$) such that
\begin{equation}\label{3.7}
\zeta_k\to\zeta\quad\text{almost everywhere in $\varOmega$ as $k\to\infty$}.
\end{equation}

Due to \eqref{3.5}, the following estimate holds for all $k\in \mathbb{N}$:
\begin{multline}\label{3.8}
\int_\varOmega \varphi_k(\zeta_k)\,\zeta_k^2\, dx=
\int_\varOmega \varphi_k(\zeta_k)\,\Big(\int_0^T \alpha (\cdot,t)\, u_k(\cdot,t)\,dt\Big)^2\, dx
\\
\le \int_\varOmega \varphi_k(\zeta_k)\Big(\int_0^T \alpha^2 (\cdot,t)\,dt \int_0^T u_k^2(\cdot,t)\,dt\Big) dx
\le C_6,
\end{multline}
where $C_6=C_1\,\|\alpha\|^2_{L^\infty(\varOmega; L^2(0,T))}$ and $C_1$ is the constant from \eqref{3.5}.
As a consequence of the Fatou lemma, this estimate together with the continuity of $\varphi$ and \eqref{3.7}  implies that
\begin{equation}\label{3.9}
\int_\varOmega \varphi(\zeta)\,\zeta^2\, dx\le C_6
\end{equation}
which is a part of \eqref{3.2}. The rest of the assertions in \eqref{3.2} are simple consequence of \eqref{3.9}.

We can also prove now that $\zeta\in L^2(\varOmega)$, and namely that
\begin{equation}\label{3.10}
\int_\varOmega \zeta^2\, dx \le TC_6,
\end{equation}
where $C_6$ is the constant from estimate \eqref{3.8}. Indeed, exactly as in \eqref{3.8}, we obtain that
$$
\int_\varOmega \zeta_k^2\, dx
\le \int_\varOmega \int_0^T \alpha^2 (\cdot,t)\,dt \int_0^T u_k^2(\cdot,t)\,dt\, dx
\le T \,\|\alpha\|^2_{L^\infty(\varOmega; L^2(0,T))} \esup_{t\in[0,T]}\|u_k(\cdot,t)\|^2
\le TC_6,
$$
and \eqref{3.10} follows from the Fatou lemma. Notice that we can put the constant
from the Poin\-care inequality instead of $T$ in \eqref{3.10}, if we estimate
$\int_\varOmega  u_k^2\, dx$ by $\int_\varOmega  |\nabla u_k|^2\, dx$ and use again \eqref{3.5}.

\medskip\noindent
\textsc{Step~3.}
In this step, we prove \eqref{3.6}, which completes the proof of the weak solvability of the problem.
Let $h:\varOmega_T\to \mathbb{R}$ be an arbitrary smooth bounded function and $R=\max_{(x,t)\in\varOmega_T}|h(x,t)|$.
We will prove that
\begin{equation}\label{3.11}
\varphi_k(\zeta_k)\int_0^T u_k\,h\,dt \to \varphi(\zeta)\int_0^T u\,h\,dt\quad\text{as $k\to\infty$ in $L^1(\varOmega)$},
\end{equation}
which clearly implies \eqref{3.6} and, as a consequence, the weak solvability of problem~\eqref{1.1}--\eqref{1.3}.
To do this, we will apply the Vitali convergence theorem (see, e.g., \cite[Sec.~4.8.7]{MP}).

Exactly as for the sequence $\{\zeta_k\}$ at the end of the proof of Theorem~\ref{t2.1},
we can establish (selecting, if necessary, a subsequence) that $\int_0^T u_k\,h\,dt\to \int_0^Tu\,h\,dt$ in
$L^2(\varOmega)$ and almost everywhere in $\varOmega$ as $k\to\infty$.
Besides that, due to  \eqref{3.7},
\begin{equation}\label{3.12}
\varphi_k(\zeta_k)\to \varphi(\zeta)\quad\text{as  $k\to\infty$ almost everywhere in $\varOmega$.}
\end{equation}
 Really, let $E\subset \varOmega$ be the set such that $\zeta_k(x)\to \zeta(x)$ as $k\to \infty$
 for all $x\in E$. As it follows from \eqref{3.7}, $\mu_n(\varOmega\setminus E)=0$,
 where $\mu_n$ is the $n$-dimensional Lebesgue measure.
For an arbitrary $x\in E$, there exists $k_{x,1}\in \mathbb{N}$ such that $\varphi(\zeta(x))+1 < k_{x,1}$.
Since the function $\varphi$ is continuous, there exists $k_{x,2}\in \mathbb{N}$ such that
$|\varphi(\zeta_k(x))-\varphi(\zeta(x))|<1$ for all $k> k_{x,2}$. Therefore,
$\varphi(\zeta_k(x))< k$ and, as a consequence, $\varphi_k(\zeta_k(x))=\varphi(\zeta_k(x))$
whenever $k>\max\{k_{x,1},k_{x,2}\}$. Due to the continuity of $\varphi$,
we conclude that $\varphi_k(\zeta_k(x))=\varphi(\zeta_k(x))\to \varphi(\zeta(x))$ as  $k\to\infty$ for all $x\in E$.
Thus,
$$
\varphi_k(\zeta_k)\int_0^T u_k\,h\,dt \to \varphi(\zeta)\int_0^T u\,h\,dt\quad\text{
almost everywhere in $\varOmega$ as $k\to\infty$}.
$$
According to the Vitali theorem, in order to prove \eqref{3.11}, it is necessary to establish that the sequence
$\{\varphi_k(\zeta_k)\int_0^T u_k\,h\,dt\}$ is uniformly integrable on $\varOmega$.

Let $E$ be an arbitrary measurable subset of $\varOmega$, $E_T=E\times [0,T]$, and
$G_M^k=\{(x,t)\in E_T\; |\; |u_k(x,t)|\ge M\}$. Then
$$
\int_{G_M^k}\varphi_k(\zeta_k)|u_k\,h|\,dx\,dt\le \frac{R}{M}\int_{G_M^k}\varphi_k(\zeta_k)\,u_k^2\,dx\,dt
\le \frac{R}{M}\int_{\varOmega_T}\varphi_k(\zeta_k)\,u_k^2\,dx\,dt\le \frac{C_1R}{M},
$$
where $C_1$ is the constant from \eqref{3.5}. Besides that,
$$
\int_{E_T\setminus G_M^k}\varphi_k(\zeta_k)|u_k\,h|\,dx\,dt\le M R\int_{E_T\setminus G_M^k}\varphi_k(\zeta_k)\,dx\,dt
\le MR\,T \int_{E}\varphi_k(\zeta_k)\,dx.
$$
Therefore,
\begin{equation}\label{3.13}
\Big|\int_{E}\varphi_k(\zeta_k)\int_0^T u_k\,h\,dt\,dx\Big|
=\Big|\int_{E_T}\varphi_k(\zeta_k)\, u_k\,h\,dx\,dt\Big|
\le\frac{C_1R}{M}+MR\, T\int_{E}\varphi_k(\zeta_k)\,dx.
\end{equation}

Let us estimate $\int_{E}\varphi_k(\zeta_k)\,dx$.
If $E_N^k=\{x\in E\; |\; |\zeta_k(x)|> N\}$, where $N\in (0,\infty)$, then \eqref{3.8} implies that
$$
\int_{E_N^k}\varphi_k(\zeta_k)\,dx\le
\frac{1}{N^2}\int_{E_N^k}\varphi_k(\zeta_k)\zeta_k^2\,dx
\le \frac{C_6}{N^2}.
$$
Since the function $\varphi$ is continuous, there exists a constant $\gamma_N>0$ such that
$\varphi(\xi)\le \gamma_N$ for all $\xi\in [-N,N]$. Therefore,
$$
\int_{E\setminus E_N^k}\varphi_k(\zeta_k)\,dx\le \int_{E\setminus E_N^k}\varphi(\zeta_k)\,dx
\le \gamma_N \mu_n(E).
$$
Thus,
$$
\int_{E}\varphi_k(\zeta_k)\,dx\le \frac{C_6}{N^2}+\gamma_N \mu_n(E)
$$
and, as it follows from \eqref{3.13},
$$
\Big|\int_{E}\varphi_k(\zeta_k)\int_0^T u_k\,h\,dt\,dx\Big|\le
\frac{C_1R}{M}+\frac{MR\, T C_6}{N^2}+\gamma_N MR\, T \mu_n(E)
$$

For an arbitrary $\varepsilon>0$, we take $M$ and $N$ such that $C_1 R/M< \varepsilon/3$ and
$MR\, T C_6/N^2<\varepsilon/3$. If $\delta= \varepsilon/(3\gamma_N MR\, T)$,
then, as a consequence of the last inequality,
$$
\Big|\int_{E}\varphi_k(\zeta_k)\int_0^T u_k\,h\,dt\,dx\Big|<\varepsilon\quad\text{whenever $\mu_n(E)<\delta$},
$$
which  implies the required uniform integrability of the sequence $\{\varphi_k(\zeta_k)\int_0^T u_k\,h\,dt\}$.
Thus, \eqref{3.11} is proven.

\medskip\noindent
\textsc{Step~4.}
In this step, we prove \eqref{3.1} and \eqref{3.3}, which completes the proof of the theorem.
The fact that $\varphi(\zeta) u^2\in L^1(\varOmega_T)$ is a consequence of \eqref{3.1}.
In order to establish \eqref{3.1}, we have to pass to the limit as $k\to\infty$ in \eqref{3.5}.
The first two terms on the left-hand side of \eqref{3.1} are obtained in the usual way due to the weak lower
semicontinuity of the norms. Let us obtain the third term.

Using the H\"older inequality, we get:
\begin{multline*}
\int_\varOmega \Big|\nabla\int_0^T u_k^2\,dt\Big|\, dx
\le 2\int_{\varOmega_T} |u_k|\,|\nabla u_k|\,dx\,dt
\\
\le 2 \Big(\int_{\varOmega_T} |u_k|^2\,dx\,dt\Big)^{1/2}\,\Big(\int_{\varOmega_T} |\nabla u_k|^2\,dx\,dt\Big)^{1/2}
\le 2T^{1/2}C_1,
\end{multline*}
where $C_1$ is the constant from \eqref{3.5}.
This estimate implies that there exists a function $v\in L^1(\varOmega)$ such that (up to a subsequence)
$$
\int_0^T u_k^2\, dt\to v\quad\text{as $k\to \infty$ in $L^1(\varOmega)$ and almost everywhere in $\varOmega$.}
$$
Therefore, due to \eqref{3.12}, $\varphi(\zeta_k)\int_0^T u_k^2\,dt\to \varphi(\zeta)\,v$
almost everywhere in $\varOmega$. As it follows from \eqref{3.5} and the Fatou lemma,
$\varphi(\zeta)\,v\in L^1(\varOmega)$ and
$$
\int_\varOmega \varphi(\zeta)\,v\,dx\le \liminf_{k\to\infty}
\int_0^T\int_\varOmega\varphi_k(\zeta_k)\, u_k^2\, dx\,dt.
$$
Thus, in order to establish \eqref{3.1}, it suffices to show that
$$
\int_0^T\int_\varOmega\varphi(\zeta)\, u^2\, dx\,dt\le \int_\varOmega \varphi(\zeta)\,v\,dx.
$$
Since $\varphi(\zeta)\ge 0$, it is enough to prove that
$$
\int_0^T u^2\,dt\le v\quad\text{almost everywhere in $\varOmega$.}
$$

Let $G$ be an arbitrary measurable subset of $\varOmega$ and $G_T=G\times [0,T]$.
Since $u_k$ converges to $u$ weakly in $L_2(G_T)$ and the norm in $L^2(G_T)$ is weakly lower semicontinuous,
$$
\int_G \int_0^T u^2\,dt\,dx\le \liminf_{k\to\infty} \int_G \int_0^T u_k^2\,dt\,dx=
\int_G v\,dx.
$$
Due to the arbitrariness of $G$, we have the required relation.

Finally, the assertion that $\varphi(\zeta) u\in L^1(\varOmega_T)$ follows from the H\"older inequality and the facts
that $\varphi(\zeta)\in L^1(\varOmega)$ and $\varphi(\zeta) u^2\in L^1(\varOmega_T)$.

The theorem is proven.
\hfill$\square$

Notice that Remark~\ref{t2.2} is also true for Theorem~\ref{t3.1}, since we used only the energy estimate
in its proof.

\bigskip\noindent
\textbf{Acknowledgement} \\
This work was supported by the Russian Science Foundation (Grant No. 19-11-00069).

\end{document}